\documentclass[12]{article}

\usepackage{amssymb}
\usepackage{color, amsmath,amssymb, amsfonts, amstext,amsthm, latexsym}

\usepackage{graphicx}
\usepackage{longtable}

\renewcommand{\phi}{\varphi}

\newcommand{\la}{{\lambda}}
\newtheorem{theorem}{Theorem}

\newtheorem{remark}{Remark}

\def\EE{{\mathbb{E}}}

\begin{document}

\title{Higher-order derivative of intersection local time for two independent
fractional Brownian motions
\thanks{Jingjun Guo  acknowledges  the support of National
Natural Science Foundation of China \#71561017, the Youth Academic Talent Plan of Lanzhou University of Finance and Economics.  Yaozhong Hu is partially supported by a grant from the Simons Foundation \#209206
and by a General Research Fund of University of Kansas.
Yanping Xiao acknowledges  the support of Northwest University for Nationalities  \#xbmuyjrc201701.}}

\author{Jingjun Guo$^{a,}\thanks{Corresponding author.
gjjemail@126.com}$, Yaozhong Hu$^{b}$, ~Yanping Xiao$^{c}$ \\
{\it\small a.School of Statistics, Lanzhou University of Finance and Economics,} \\
{\it\small Lanzhou,GS 730020,P.R.China }\\
{\it\small b.Department of Mathematics, University of Kansas,}\\
{\it\small Lawrence, KS 66045, USA}\\
{\it\small c.School of Mathematics and Computer Science, }\\
{\it\small Northwest University for Nationalities, Lanzhou,GS 730000, P.R.China}}
\date{}

\maketitle

\begin{abstract}
In this article, we obtain sharp conditions for the
existence of the high order  derivatives ($k$-th order)
of intersection local time
$ \widehat{\alpha}^{(k)}(0)$ of
 two independent d-dimensional fractional Brownian motions $B^{H_1}_t$ and
 $\widetilde{B}^{H_2}_s$ with Hurst parameters $H_1$ and
 $H_2$, respectively.  We also study their exponential integrability.



\bigskip

\end{abstract}

\section{Introduction and main result}\label{sec-1}
Intersection local time or self-intersection local time when the two processes
%fractional Brownian motions
are the same are important subjects  in probability theory
and their  derivatives have received  much attention  recently.
Jung et al.  \cite{jung1} and \cite{jung2}   discussed
Tanaka formula and occupation-time formula for derivative self-intersection local time
of  fractional Brownian motions.
On the other hand,  several  authors paid  attention to the  renormalized self-intersection local time of  fractional Brownian motions, see  e.g.,
Hu et al.   \cite{hunualart} and \cite{hunualartsong}.

Motivated by \cite{jung1} and \cite{hu17}, higher-order derivative    of intersection local time for two independent
fractional Brownian motions is studied in this paper.

To state our main result we
let $B^{H_1}=\{B^{H_1}_t, t\geq 0\}$ and $\widetilde{B}^{H_2}=\{\widetilde{B}^{H_2}_t, t
\geq 0
\}$ be two independent $d$-dimensional
fractional Brownian motions of Hurst parameters $H_1, H_2\in (0, 1)$ respectively. This means that
$B^{H_1}$ and $\widetilde{B}^{H_2}$
are  independent   centered Gaussian processes with covariance
$$\mathbb{E}[B^{H_{1 } }_s   B^{H_{1} }_t]=\frac{1}{2}(s^{2H_{1}}+t^{2H_{1}}-{\mid
s-t
\mid}^{2H_1})\,.  $$
(similar identity for $\tilde B$).

In this paper we are concerned with
the  derivatives of  intersection local time of
$B^{H_1}$ and $\widetilde{B}^{H_2}$,    defined by
\begin{equation*}
\hat{\alpha}^{(k)}(x):= \frac{\partial ^k}{\partial x_1^{k_1}\cdots\partial x_d^{k_d}}
\int_0^T\int_0^T\delta(B^{H_1}_t-\widetilde{B}^{H_2}_s+x)dtds,
\end{equation*}
where $k=(k_1, \cdots, k_d)$ is a multi-index with all
$k_i$ being nonnegative integers and $\delta$ is the Dirac delta
function of $d$-variable.  In particular, we are exclusively consider
the case when $x=0$ in this work. Namely, we are studying
\begin{equation}
\hat{\alpha}^{(k)}(0):=
\int_0^T\int_0^T\delta^{(k)} (B^{H_1}_t-\widetilde{B}^{H_2}_s )dtds,
\end{equation}
where $\delta^{(k)}(x)=\frac{\partial ^k}{\partial x_1^{k_1}\cdots\partial x_d^{k_d}}\delta(x)$ is k-th  order partial
derivative of the Dirac delta function.  Since $\delta(x)=0$ when $x\not =0$
the intersection local time $\hat{\alpha} (0)$   (when $k=0$)
  measures the frequency
that processes $B^{H_1}$ and $\widetilde{B}^{H_2}$ intersect each other.

Since the Dirac delta function $\delta$  is a generalized function,
we need to give a meaning to  $\hat{\alpha}^{(k)}(0)$.  To this end, we approximate the Dirac delta function $\delta$ by
\begin{equation}
f_{\varepsilon}(x):=\frac{1}{(2\pi\varepsilon)^{\frac{d}{2}}}
e^{-\frac{\mid x\mid^2}{2\varepsilon}}
=\frac{1}{(2\pi)^d}\int_{\mathbb{R}^d}e^{ipx}e^{-\frac{\varepsilon \mid p\mid^2}{2}}dp\,,
\end{equation}
and throughout this paper, we use $px=\sum_{j=1}^d p_jx_j$
and $|p|^2=\sum_{j=1}^d p_j^2$.
Thus we approximate $\delta^{(k)}$ by
\begin{equation}
\begin{aligned}
f_{\varepsilon}^{(k)}(x):=\frac{\partial ^k}{\partial x_1^{k_1}\cdots\partial x_d^{k_d}} f_{\varepsilon}(x)
=\frac{i^k}{(2\pi)^d}\int_{\mathbb{R}^d}  p_1^{k_1}\cdots p_d^{k_d}
e^{ipx}e^{-\frac{\varepsilon
\mid p\mid^2}{2}}dp.
\end{aligned}
\end{equation}

We say that $\hat{\alpha}^{(k)} (0)$ exists (in $L^2$) if
\begin{equation}
\hat{\alpha}^{(k)}_{\varepsilon}(0):= \int_0^T\int_0^Tf^{(k)}_{\varepsilon}(B^{H_1}_t-
\widetilde{B}^{H_2}_s)dtds
\end{equation}
converges to a random variable (denoted by $\hat{\alpha}^{(k)} (0)$)
in $L^2$ when $\varepsilon\downarrow 0$.

Here is the main result of this work.

\begin{theorem}
 Let $B^{H_1}$ and $\widetilde{B}^{H_2}$   be two independent  d-dimensional fractional
 Brownian motions of Hurst parameter $H_1$ and $H_2$, respectively.
 \begin{enumerate}
 \item[(i)] Assume    $k =(k_1, \cdots, k_d)$
 is an index of  nonnegative  integers (meaning that $k_1, \cdots, k_d$
 are nonnegative integers)    satisfying
 \begin{equation}
 \frac{H_1H_2}{H_1+H_2}
(|k|+d)< 1\,,\label{e.cond_H}
\end{equation}
    where $|k|=k_1+\cdots+k_d$.   Then, the k-th order derivative intersection
   local time $\hat{\alpha}^{(k)}(0)$ exists in $L^p(\Omega)$
   for any $p\in [1, \infty)$.
   \item[(ii)]
   Assume condition \eqref{e.cond_H} is satisfied. There
 is a strictly positive constant $C_{d,k,T}\in (0, \infty)$ such that
\begin{equation*}
\begin{aligned}
\mathbb{E}
\left[\exp\left\{C_{d,k,T}\left| \widehat{\alpha}^{(k)}(0)\right| ^{\beta}\right\}
\right]<\infty,
\end{aligned}
\end{equation*}
 where $\beta=\frac{H_1+H_2}{2dH_1H_2}$.
 \item[(iii)] If $\hat{\alpha}^{(k)}(0)\in L^1(\Omega)$,
 where $k=(0, \cdots, 0,k_i, 0,\cdots,0)$ with $k_i$
 being even integer,   then condition
 \eqref{e.cond_H} must be satisfied.
 \end{enumerate}
\end{theorem}

 \begin{remark}
 \begin{enumerate}
 \item[(i)] When $k=0$, we have that $\widehat{\alpha}^{(0)}(0)$ is in
 $L^p(\Omega)$ for any $p\in [1, \infty)$ if  $\frac{H_1H_2}{H_1+H_2}
d<1$.  In the special case  $H_1=H_2=H$,  this condition becomes $Hd<2$, which is the condition obtained in   Nualart et al. \cite{nualartortiz}.
 \item[(ii)]  When $H_1=H_2=\frac12$,  we have the exponential integrability
 exponent   $\beta=2/d$, which  implies an earlier result (\cite[Theorem 9.4]{hu17}).
 \item[(iii)] We also show that condition \eqref{e.cond_H} is necessary
 in some sense. This is also first time.
 \end{enumerate}
 \end{remark}

\section{Proof of the theorem}
\noindent{\bf Proof of Parts  (i) and (ii)}.
 This section is devoted to the proof of the theorem.
 We shall first find a good bound for $\EE\left| \widehat{\alpha}^{(k)}(0)\right| ^n$
 which  gives  a proof for (i) and (ii) simultaneously.
 We introduce the following notations.
\begin{eqnarray*}
&& p_j=(p_{1j},\ldots,p_{dj})\,,\qquad  p_j^{k }=(p_{1j}^{k_1 },\ldots,p_{dj}^{k_d})\,,\qquad j=1, 2, \cdots, n\,; \\
&& p=(p_1, \dots, p_n)\,,  \qquad dp =\prod_{i=1}^d\prod_{j=1}^n dp_{ij}\,.
\end{eqnarray*}
 We also denote
$ds=ds_1\cdots ds_n$ and
$dt=dt_1\cdots dt_n$.

Fix an integer $n\geq 1$. Denote  $T_n=\{0<t,s <T\}^n$.
We have
\begin{eqnarray*}
\begin{aligned}
\mathbb{E}&\left[\left|\widehat{\alpha}^{(k)}_{\varepsilon}(0)\right|^n
\right]\\
&\le \frac{1}{(2\pi)^{nd}}\int_{T_n}\int_{\mathbb{R}^{nd}} \bigg| \mathbb{E}[\exp\{ip_1(B^{H_1}_{s_1}-\widetilde{B}^{H_2}_{t_1})+\cdots
\\
&~~+ip_n(B^{H_1}_{s_n}-\widetilde{B}^{H_2}_{t_n})\}]\bigg| \exp\{-\frac{\varepsilon}{2}\sum_{j=1}
^n\mid p_j\mid^2\}\prod_{j=1}^n
\mid p_j ^k\mid  dp dtds\\
&=\frac{1}{(2\pi)^{nd}}\int_{T_n}\int_{\mathbb{R}^{nd}}  \exp\{-\frac{1}{2}\mathbb{E}\big[\sum_{j=1}^np_j(B^{H_1}_{s_j}-
\widetilde{B}^{H_2}_{t_j})
\big]^2\}
\\
&~~\times\exp\{-\frac{\varepsilon}{2}\sum_{j=1}
^n\mid p_j\mid^2\} \prod_{j=1}^n
\mid p_j ^k\mid  dpdtds\\
&\le \frac{1} {(2\pi)^{nd}}\int_{T_n}\int_{\mathbb{R}^{nd}}    \prod_{i=1}^d  \left(\prod_{j=1}^n
\mid p_{ij}^{k_i}\mid  \right)     \exp\{-\frac{1}{2}\mathbb{E}[p_{i1}B^{H_1,i}_{s_1}+\cdots+p_{in}B^{H_1,i}_{s_n}]^2\\
&~~-\frac{1}{2}\mathbb{E}[p_{i1}B^{H_2,i}_{t_1}+\cdots+p_{in}B^{H_2,i}_{t_n}]^2\}dpdtds \,.
\end{aligned}
\end{eqnarray*}
  The expectations in the
above exponent  can be computed by
\begin{eqnarray*}
&&\mathbb{E}[p_{i1}B^{H_1,i}_{s_1}+\cdots+p_{in}B^{H_1,i}_{s_n}]^2=
(p_{i1}, \cdots, p_{in}) Q_1(p_{i1}, \cdots, p_{in})^T\,,\\
&&
\mathbb{E}[p_{i1}\tilde B^{H_2,i}_{s_1}+\cdots+p_{in}
\tilde B^{H_2,i}_{s_n}]^2=(p_{i1}, \cdots, p_{in}) Q_2 (p_{i1}, \cdots, p_{in})^T\,,
\end{eqnarray*}
where
\[
Q_1=\EE \left(B^{H_1, i}_j B^{H_1, i}_k\right)_{1\le j,k\le n}
\quad {\rm and}\quad  Q_2=\EE \left(\tilde B^{H_2, i}_j\tilde
B^{H_2, i}_k\right)_{1\le  j,k\le n}
\]
 denote respectively covariance matrices  of n-dimensional random vectors
$(B^{H_1,i}_{s_1},...,B^{H_1,i}_{s_n})$ and that of
$(\widetilde{B}^{H_2,i}_{t_1},...,
\widetilde{B}^{H_2,i}_{t_n})$.  Thus we have
\begin{eqnarray}
\mathbb{E}\left[
\left|\widehat{\alpha}^{(k)}_{\varepsilon}(0)\right|^n
\right]\le \frac{1}{(2\pi)^{nd}}  \int_{T_n}  \prod_{i=1}^d  I_i(t,s)  dtds\,,
\label{e.boundedbyI}
\end{eqnarray}
where
\[
I_i(t,s):= \int_{\mathbb{R}^{n}}
\mid x^{k_i} \mid \exp\{-\frac{1}{2}x^T(Q_1+Q_2)x\} dx  \,.
\]
Here we recall $x=(x_1, \cdots, x_n)$  and
$x^k_i=x_1^{k_i} \cdots x_n^{k_i}$.

For each fixed $i$ let us compute integral $I_i(t,s)$ first.  Denote  $B=Q_1+Q_2$.
Then $B$ is a strictly positive definite matrix and hence $\sqrt B$ exists.  Making substitution
$\xi=\sqrt{B}x$.  Then
\begin{eqnarray*}
I_i(t,s)=\int_{\mathbb{R}^{n }}\prod_{j=1}^n
\mid (B^{-\frac{1}{2}}\xi)_j\mid^{k_i} \exp\{-\frac{1}{2}\mid \xi\mid^2\}\mathrm{det}(B)
^{-\frac{1}{2}}d\xi.\\
\end{eqnarray*}
To obtain a nice bound for the above
integral, let us first diagonalize   $B$:
\begin{equation*}
B =Q \Lambda Q^{-1} \,,
\end{equation*}
where $\Lambda=$diag$\{\lambda_1 ,...,\lambda_n \}$ is a strictly
positive diagonal matrix with $\la_1\le \la_2\le\cdots\le \la_d$
and    $Q =(q_{ij})_{1\le i,j\le d} $ is  an orthogonal matrix.  Hence, we have
$\det(B)=\la_1\cdots\la_d$.
Denote
\begin{equation*}
\eta=
\begin{pmatrix}
\eta_1, \eta_2, \cdots,
\eta_n
\end{pmatrix}^T=Q^{-1}\xi,
\end{equation*}
Hence,
\begin{eqnarray*}
B^{-\frac{1}{2}}\xi
&=&Q\Lambda^{-1/2}Q^{-1} \xi =
Q\Lambda^{-1/2}\eta\\
&=&Q\begin{pmatrix}
\lambda_1^{-\frac{1}{2}}\eta_1\\\lambda_2^{-\frac{1}{2}}\eta_2\\\vdots\\\lambda_n^{-\frac{1}{2}}
\eta_n
\end{pmatrix}
=
\begin{pmatrix}
q_{1,1}&q_{1,2}&\cdots&q _{1,n} \\
q_{2,1}&q_{2,2}&\cdots&q_{2,n}\\
\vdots&\vdots&\cdots&\vdots\\
q_{n,1}&q _{n,2}&\cdots&q_{n,n}\\
\end{pmatrix}
\begin{pmatrix}
\lambda_1^{-\frac{1}{2}}\eta_1\\\lambda_2^{-\frac{1}{2}}\eta_2\\\vdots\\\lambda_n^{-\frac{1}{2}}
\eta_n.
\end{pmatrix}.
\end{eqnarray*}

Therefore, we have
\begin{eqnarray*}
\mid(B^{-\frac{1}{2}}\xi)_j\mid
&=&\mid\sum_{k=1}^nq_{jk}\lambda_{k}^{-\frac{1}{2}}\eta_{k}\mid
\leq
\lambda_1
^{-\frac{1}{2}}  \sum_{k=1}^n\mid q_{jk} \eta_{k}\mid \\
&\leq & \lambda_1
^{-\frac{1}{2}}
\left(\sum_{k=1}^nq_{jk}^2\right)^{\frac{1}{2}}
\left(\sum_{k=1}^n\eta_{k}^2\right)^{\frac{1}{2}}
 \leq  \lambda_1
^{-\frac{1}{2}} \mid\eta\mid_2
 =\lambda_1
^{-\frac{1}{2}}\mid\xi\mid_2.\\
\end{eqnarray*}
 Since both $Q_1$ and $Q_2$ are positive definite, we see that
 \[
 \la_1\ge \la_1(Q_1) \,,\qquad {\rm and}\qquad  \la_1\ge \la_1(Q_2)\,,
 \]
 where $  \la_1(Q_i)$ is the smallest eigenvalue of
 $Q_i$, $i=1, 2$.   This means that
\begin{equation*}
 \la_1 \ge \la_1(Q_1)^\rho \la_1(Q_2)^{1-\rho}\quad \hbox{for any }\
 \ \rho\in [0, 1]\,.
 \end{equation*}
 This implies
 \begin{equation*}
\mid(B^{-\frac{1}{2}}\xi)_j\mid
\le \la_1(Q_1)^{-\frac12\rho} \la_1(Q_2)^{-\frac12(1-\rho)}
\mid\xi\mid_2.\label{e.Bxi_bound}
\end{equation*}
Consequently, we have
\begin{eqnarray}
I_i(t,s)&=& \mathrm{det}(B)
^{-\frac{1}{2}} \la_1(Q_1)^{-\frac12\rho k_i } \la_1(Q_2)^{-\frac12(1-\rho)k_i } \nonumber\\
&&\qquad   \int_{\mathbb{R}^{n }}
\mid\xi\mid_2 ^{k_i}  \exp\{-\frac{1}{2}\mid \xi\mid^2\} d\xi\,,
\label{e.I_first_bound}
\end{eqnarray}
for any $\rho\in [0,1]$.

Now we are going to find a lower bound for $\la_1(Q_1)$ ($\la_1(Q_2)$ can be dealt with the same way. We only need to replace $s$ by $t$).   Without loss of generality we can assume $0\le s_1<s_2<\cdots<s_n\le T$.  From the definition of $Q_1$ we have for any vector $u=(u_1, \cdots, u_d)^T$,
\begin{equation*}
\begin{aligned}
&
u^TQ_1u
~=\mathrm{Var}\big(u_{1}B^{H_1}_{s_{1}}+u_{2}B^{H_1}_{s_{2}}
+...+
u_{n}B^{H_1}_{s_{n}}\big)\\
&~=\mathrm{Var}\big((u_{1}+...+u_{n})B^{H_1}_{s_{1}}
+(u_{2}+...+u_{n})(B^{H_1}_{s_{2}}-B^{H_1}_{s_{1}})\\
&~~+...+
(u_{n-1}+u_{n})(B^{H_1}_{s_{n-1}}-B^{H_1}_{s_{n-2}})+
u_{n}(B^{H_1}_{s_{n}}-B^{H_1}_{s_{n-1}})\big)\\
\end{aligned}
\end{equation*}
Now  using Lemma 8.1 of  \cite{nondet} and Lemma A.1 of \cite{hunualartsong}, 
we see that $R$ appeared in Lemma 8.1 satisfies $R\ge \prod_{i=1}^n \sigma_i^2 $
($m$ in Lemma 8.1 of \cite{nondet} is replaced by $n$).  Thus
we have
\begin{equation*}\begin{aligned}
u^TQ_1u&~\geq  \big((u_{1}+...+u_{n})^2s_{1}^{2H_1}
+(u_{2}+...+u_{n})^2(s_{2}-s_{1})^{2H_1}\\
&~~+...+
(u_{n-1}+u_{n})^2(s_{n-1}-s_{n-2})^{2H_1}+
u_{n}^2(s_{n}-s_{n-1})^{2H_1}\big)\\
&~\geq \min\{s_{1}^{2H_1},(s_{2}-s_{1})^{2H_1},...,
(s_{n}-s_{n-1})^{2H_1}\}\\
&~~\big[(u_{1}+...+u_{n})^2+(u_{2}+...+u_{n})^2
+...+(u_{n-1}+u_{n})^2+u_{n}^2\big]\,.
\end{aligned}
\end{equation*}

Consider the function
\begin{eqnarray*}
f(u_1, \cdots,  u_n)
&=&(u_{1}+\cdots +u_{n})^2+(u_{2}+\cdots +u_{n})^2
+\cdots +(u_{n-1}+u_{n})^2+u_{n}^2 \\
&=&(u_1, \cdots, u_n) G(u_1, \cdots, u_n)^T\,,
\end{eqnarray*}
where
\[
G=\left(\begin{matrix}1&1&1&\cdots&1\\ 0&1&1&\cdots&1\\
\vdots&\vdots&\vdots& \cdots&\vdots\\
0&0&0&\cdots&1
\end{matrix}\right)\,.
\]
 It is easy to see that the matrix $G^TG$ has a minimum eigenvalue independent of $n$. Thus
this function $f$ attains its minimum value $f_{\rm min} $ independent of $n$
 on the sphere
$u_1^2+\cdots+u_n^2=1$.   It is also easy to see that $f_{\rm min}>0$.

As a consequence we have
\begin{eqnarray}
&&\la_1(Q_1)
=\inf_{|u|=1} u^TQ_1u\nonumber\\
&&\qquad \ge    \min\{s_{1}^{2H_1},(s_{2}-s_{1})^{2H_1},...,
(s_{n}-s_{n-1})^{2H_1}\}\inf _{|u|=1} f(u_1, \cdots, u_n)\nonumber\\
&&\qquad \ge   f_{\rm min} \min\{s_{1}^{2H_1},(s_{2}-s_{1})^{2H_1},...,
(s_{n}-s_{n-1})^{2H_1}\}\nonumber\\
&&\qquad \ge   K \min\{s_{1}^{2H_1},(s_{2}-s_{1})^{2H_1},...,
(s_{n}-s_{n-1})^{2H_1}\}\,. \label{e.la_1_q1}
\end{eqnarray}
In a similar way we have
\begin{eqnarray}
\la_1(Q_2)
%&=&\inf_{|u|=1} u^TQ_1u\\
%&\ge&  K\min\{s_{1}^{2H_1},(s_{2}-s_{1})^{2H_1},...,
%(s_{n}-s_{n-1})^{2H_1}\}\inf _{|u|=1} f(u_1, \cdots, u_n)\\
%&\ge& Kf_{\rm min} \min\{s_{1}^{2H_1},(s_{2}-s_{1})^{2H_1},...,
%(s_{n}-s_{n-1})^{2H_1}\}\\
&\ge&  K \min\{t_{1}^{2H_2},(t_{2}-t_{1})^{2H_2},...,
(t_{n}-t_{n-1})^{2H_2}\}\,. \label{e.la_1_q2}
\end{eqnarray}
%As a result, it exists
%\begin{equation*}
%\begin{aligned}
%\min\{\lambda_1,...,\lambda_n\}\geq K\min_{i,j=1,...,n}\{( s_{\sigma(i)}-s_{i-1})^{2H_1}+(
%t_{j}-t_{j-1})^{2H_2}\},\\
%\end{aligned}
%\end{equation*}
%where $s_{\sigma(0)}\equiv0$, $t_{\sigma(0)}\equiv0$ and K is also some constant.
The integral in \eqref{e.I_first_bound}   can be bounded as
\begin{eqnarray}
I_2&:=&\int_{\mathbb{R}^{n }}
 \mid\xi \mid  ^{ k_i }
\exp\{-\frac{1}{2} \mid\xi \mid^2 \}
d\xi \nonumber\\
&\leq&n^{\frac{k_i}{2}}\int_{\mathbb{R}^{nd}}\max_{1\le j\le n}  \mid \xi_j\mid^{k_i}  \exp\{-\frac{1}{2}
 \mid\xi \mid^2 \}
d\xi \nonumber\\
&\leq&n^{\frac{k_i}{2}}\int_{\mathbb{R}^{n }}\sum_{j=1}^n \mid \xi_j\mid^{k_j}  \exp\{-\frac{1}{2} \mid\xi \mid^2 \}
d\xi \nonumber\\
&\leq&n^{\frac{k_i}{2}+1}\int_{\mathbb{R}^{n }} \mid \xi_{1}\mid^{k_i}  \exp\{-\frac{1}{2}
 \mid\xi \mid^2 \}
d\xi\nonumber\\
& \leq&  n^{\frac{k_i}{2}+1}C^n \le C^n\,. \label{e.I_2}
\end{eqnarray}
%\begin{eqnarray*}
%I_2&:=&\int_{\mathbb{R}^{n }}
%(\mid\xi_1\mid^2+\cdots+\mid\xi_n\mid^2)^{\frac{k_i}{2}}
%\exp\{-\frac{1}{2}(\mid\xi_1\mid^2+\cdots+\mid\xi_n\mid^2)\}
%d\xi_1\cdots d\xi_n\\
%&\leq&n^{\frac{k}{2}}\int_{\mathbb{R}^{nd}}\max\{\mid \xi_j\mid^k\}\exp\{-\frac{1}{2}
%(\mid\xi_1\mid^2+\cdots+\mid\xi_n\mid^2)\}
%d\xi_1\cdots d\xi_n\\
%&\leq&n^{\frac{k}{2}}\int_{\mathbb{R}^{nd}}\sum_{j=1}^n\{\mid \xi_1\mid^k+\cdots+\mid
%\xi_n\mid^k\}\exp\{-\frac{1}{2}(\mid\xi_1\mid^2+\cdots+\mid\xi_n\mid^2)\}
%d\xi_1\cdots d\xi_n\\
%&\leq&n^{\frac{k}{2}+1}\int_{\mathbb{R}^{nd}}\{\mid \xi_{j_0}\mid^k\}\exp\{-\frac{1}{2}
%(\mid\xi_1\mid^2+\cdots+\mid\xi_n\mid^2)\}
%d\xi_1\cdots d\xi_n\\
%&\leq &n^{\frac{n}{2}+1}C^n\\
%&\leq &C^n,\\
%\end{eqnarray*}
Substitute \eqref{e.la_1_q1}-\eqref{e.I_2}  into \eqref{e.I_first_bound}
%where $\mid \xi_{j_0}\mid$ denotes the maximum of $\mid\xi_{1}\mid^k,...,\mid\xi_{n}\mid^k$.
%Combining this with the fact
%
%\begin{equation*}
%\int_{\mathbb{R}}x^ne^{-\frac{1}{2}x^2}dx=
%  \left\{
%   \begin{aligned}
%   &\sqrt{2\pi}(n-1)!!, if~~ n ~~is ~~even,  \\
%   &0, ~~~~~~~~~~~~~~~~if ~~n ~~is ~~odd,\\
%   \end{aligned}
%   \right.
%  \end{equation*}
%for all $n\geq0$, it shows that $\Delta_4$ is finite.
%
%
%Hence,
we obtain
\begin{eqnarray*}
I_i(t,s)&\le & C^n \mathrm{det}(B)
^{-\frac{1}{2}}  \min_{j=1,...,n}(s_{j}-s_{ j-1 })^{- \rho H_1k_i} \nonumber\\
&&\qquad \min_{j=1,...,n}(t_{j}-t_{ j-1 })^{- (1-\rho)  H_2k_i} \,.
\end{eqnarray*}
%Choose $\rho=\frac{H_2}{H_1+H_2}$, we see that
%\begin{eqnarray}
%I_i(t,s)&\le & C^n \mathrm{det}(B)
%^{-\frac{1}{2}}  \left[\min_{j=1,...,n}(s_{j}-s_{ j-1 }) \min_{j=1,...,n}(t_{j}-t_{ j-1 })\right]^{- \frac{H_1  H_2} {H_1+H_2}}\,. \nonumber\\
%\label{e.I_second_bound}
%\end{eqnarray}
%
%
%\begin{equation*}
%\begin{aligned}
%\Delta_3
%\leq C^n \mathrm{det}(B)^{-\frac{d}{2}}\left[\min_{i=1,...,n}\{ s_{\sigma(i)}-s_{i-1}\}
%\min_{j=1,...,n}\{t_{j}-t_{j-1}\}\right]^{-\frac{2H_1H_2dk}{H_1+H_2}}.
%\end{aligned}
%\end{equation*}
%
%Therefore,
%\begin{eqnarray}
%I_i(t,s)&\leq& K\mathrm{det}(B)^{-\frac{d}{2}}
%\min_{i,j=1,...,n}\{ (s_{\sigma(i)}-s_{(i-1)})^{2H_1}
%+(t_{j}-t_{j-1})^{2H_2}\}^{-dk}\\
%&&~~\times\int_{\mathbb{R}^{nd}}
%(\mid\xi_1\mid^2+\cdots+\mid\xi_n\mid^2)^{\frac{k}{2}}
%\exp\{-\frac{1}{2}(\mid\xi_1\mid^2+\cdots+\mid\xi_n\mid^2)\}
%d\xi_1\cdots d\xi_n,\\
%\end{eqnarray}
%when k is even. We estimate respectively three parts in the right side of (13).
%
%Indeed, it is easily to get
%\begin{equation*}
%\begin{aligned}
%&\min_{i,j=1,...,n}\{(s_{\sigma(i)}-s_{i-1})^{2H_1}+( t_{j}-t_{j-1})^{2H_2}\}\\
%&~~\geq\min_{i,j=1,...,n}\{(s_{\sigma(i)}-s_{i-1})^{2H_1},( t_{j}-t_{j-1})^{2H_2}\}\\
%&~~=\min_{i=1,...,n}\{(s_{\sigma(i)}-s_{i-1})^{2\gamma H_1}\}\min_{j=1,...,n}\{( t_{j}-t_{j-1})^{2(1-\gamma)H_2}\},\\
%\end{aligned}
%\end{equation*}
%where we take $\gamma=\frac{H_2}{H_1+H_2}.$
%

Next we  obtain a lower bound for    $\mathrm{det}(B) $.   According to
   \cite[Lemma 9.4]{hu17}
$$\mathrm{det}(Q_1+Q_2)\geq \mathrm{det}(Q_1)^{\gamma}\mathrm{det}(Q_2)^{1-\gamma},$$
 for  any two symmetric positive definite matrices  $Q_1$ and $Q_2$
 and for any $\gamma\in [0, 1]$.
%we have
%$$\mathrm{det}(B)^{-\frac{d}{2}}= \mathrm{det}(Q_1+Q_2)^{-\frac{d}{2}}\leq \mathrm{det}(Q_1)^
%{-\frac{d}{4}}\mathrm{det}(Q_2)^{-\frac{d}{4}}.$$
%So, the following estimation is obtained
%\begin{equation*}
%\begin{aligned}
%\Delta_3\leq& C^n\left[\min_{i=1,...,n}\{s_{\sigma(i)}-s_{i-1}\}\min_{j=1,...,n}\{ t_{j}-
%t_{j-1}\}\right]^{-\frac{2H_1H_2dk}{H_1+H_2}}\\ &~~\times\mathrm{det}(Q_1)^
%{-\frac{d}{4}}\mathrm{det}(Q_2)^{-\frac{d}{4}}.
%\end{aligned}
%\end{equation*}
Now it is well-known that (see
also the usages in  \cite{hu17}, \cite{hunualart} and \cite{hunualartsong}).
%, we introduce
%$\mathfrak{F}_i$ to denote $\sigma$-algebra formation of $\{B_{s_{1}}^{H_1},B_{s_{2}}^{H_1},...,
%B_{s_{i-1}}^{H_1}\}$ and obtain
\begin{equation*}
\begin{aligned}
\mathrm{det}(Q_1)
%\mathrm{Var}(B^{H_1}_{s_{1}})\mathrm{Var}(B^{H_1}_{s_{2}}\mid B^{H_1}_{s_{1}})
%\cdots
%\mathrm{Var}(B^{H_1}_{s_{n}}\mid B^{H_1}_{s_{1}}, B^{H_1}_{s_{2}},...,B^{H_1}_{s_{n-1}})\\
%\geq& \mathrm{Var}(B^{H_1}_{s_{1}})\mathrm{Var}(B^{H_1}_{s_{2}}\mid \mathfrak{F}_{s_{1}})\cdots
%\mathrm{Var}(B^{H_1}_{s_{n}}\mid  \mathfrak{F}_{s_{n-1}})\\
\geq& C^ns_{1}^{2H_1}(s_{2}-s_{1})^{2H_1}\cdots (s_{n}-s_{n-1})^{2H_1}.\\
\end{aligned}
\end{equation*}
%Similar technique to $Q_2$, using  $\mathfrak{F}_i'$ to denote the  $\sigma$-algebra of
%$\{B_{t_{1}}^{H_2},B_{t_{2}}^{H_2},...,B_{t_{i-1}}^{H_2}\}$, we have
and
 \begin{equation*}
\begin{aligned}
\mathrm{det}(Q_2)
&\geq C^nt_{1}^{2H_2}(t_{2}-t_{1})^{2H_2}\cdots (t_{n}-t_{ n-1 })^{2H_2}.
\end{aligned}
\end{equation*}
As a consequence, we  have
\begin{eqnarray*}
I_i(t,s)&\le & C^n  \min_{j=1,...,n}(s_{j}-s_{ j-1 })^{- \rho H_1k_i}   \min_{j=1,...,n}(t_{j}-t_{ j-1 })^{- (1-\rho)  H_2k_i} \\
&&\quad \left[ s_{1} (s_{2}-s_{1}) \cdots (s_{n}-s_{n-1}) \right]^{-\gamma H_1}
\left[ t_{1} (t_{2}-t_{1}) \cdots (t_{n}-t_{n-1}) \right]^{-(1-\gamma)H_2 }.
\end{eqnarray*}

Thus,
\begin{eqnarray*}
&&\mathbb{E} \left[\left| \widehat{\alpha}^{(k)}_{\varepsilon}(0)\right|^n
\right]
\le (n!)^2 C^n  \int_{\Delta_n^2}\min_{j=1,...,n}(s_{j}-s_{ j-1 })^{- \rho H_1|k|}\\
  &&\qquad \qquad \min_{j=1,...,n}(t_{j}-t_{ j-1 })^{- (1-\rho)  H_2|k|}
  \left[ s_{1} (s_{2}-s_{1}) \cdots (s_{n}-s_{n-1}) \right]^{-\gamma H_1d}\\
&&\qquad\qquad  \left[ t_{1} (t_{2}-t_{1}) \cdots (t_{n}-t_{n-1}) \right]^{-(1-\gamma)H_2d }dtds\\
&&\qquad \le (n!)^2 C^n \sum_{i,j=1}^n  \int_{\Delta_n^2} (s_{i}-s_{ i-1 })^{- \rho H_1|k|}\\
  &&\qquad \qquad   (t_{j}-t_{ j-1 })^{- (1-\rho)  H_2|k|}
  \left[ s_{1} (s_{2}-s_{1}) \cdots (s_{n}-s_{n-1}) \right]^{-\gamma H_1d}\\
&&\qquad\qquad  \left[ t_{1} (t_{2}-t_{1}) \cdots (t_{n}-t_{n-1}) \right]^{-(1-\gamma)H_2d }dtds\,,
\end{eqnarray*}
where $\Delta_n=\left\{ 0<s_1<\cdots<s_n\le T\right\}$ denotes the
simplex in $[0, T]^n$.   We choose $\rho=\gamma=\frac{H_2}{H_1+H_2}$  to obtain \begin{eqnarray*}
&&\mathbb{E} \left[\left| \widehat{\alpha}^{(k)}_{\varepsilon}(0)\right|^n
\right]
  \le (n!)^2 C^n \sum_{i,j=1}^n   I_{3,i}I_{3,j}\,,
\end{eqnarray*}
where
\begin{eqnarray*}
I_{3,j}=   \int_{\Delta_n }    (t_{j}-t_{ j-1 })^{- \frac{H_1H_2}{H_1+H_2}
|k|}
    \left[ t_{1} (t_{2}-t_{1}) \cdots (t_{n}-t_{n-1}) \right]^{-\frac{H_1H_2}{H_1+H_2} d }dt \,,
\end{eqnarray*}
By Lemma 4.5 of \cite{HHNT}, we see that if
\[
 \frac{H_1H_2}{H_1+H_2}
(|k|+d)\le 1\,,
\]
then
\[
I_{3,j}\le \frac{ C^n T^{ \kappa_1 n-\frac{H_1H_2|k|}{H_1+H_2}  }}{\Gamma( n\kappa_1- \frac{H_1H_2}{H_1+H_2}
|k| +1)}\,,
\]
where
\[
\kappa_1=1-\frac{dH_1H_2}{H_1+H_2} \,.
\]

Substituting this bound we obtain
\begin{eqnarray*}
&&\mathbb{E} \left[\left| \widehat{\alpha}^{(k)}_{\varepsilon}(0)\right|^n
\right]
  \le n^2 (n!)^2 C^n  \frac{   T^{ 2\kappa_1 n-\frac{2H_1H_2|k|}{H_1+H_2}  }}{\Gamma^2( n\kappa_1- \frac{H_1H_2}{H_1+H_2}
|k| +1)}\\
&&\qquad \le  (n!)^2 C^n  \frac{   T^{ 2\kappa_1 n-\frac{2H_1H_2|k|}{H_1+H_2}  }}{\left(\Gamma  ( n\kappa_1- \frac{H_1H_2}{H_1+H_2}
|k| +1)\right)^2 }
\\
 &&\qquad \le  C_T (n!)^{2-2\kappa_1} C^n     T^{ 2\kappa_1 n }   \,,
\end{eqnarray*}
where $C$ is a constant independent of $T$ and $n$ and $C_T$ is a constant independent of $n$.

%
%
%According to those  equivalents above and the fact $n\leq n!$, it implies that
%\begin{equation}
%\begin{aligned}
%\mathbb{E}[(\widehat{\alpha}^{(k)}(0))^n]
%&\leq C_1T^{n(2-\frac{1}{2}dH_1-\frac{1}{2}dH_2)-\frac{4H_1H_2kd}{H_1+H_2}}(n!)^{1+\frac{1}{2}dH_1+\frac{1}{2}dH_2},\\
%\end{aligned}
%\end{equation}
%Set $\beta=\frac{2}{2+dH_1+dH_2},$ then
For any $\beta>0$, the above inequality implies
\begin{equation*}
\begin{aligned}
\mathbb{E}\left[\left|\widehat{\alpha}^{(k)}(0)\right| ^{n\beta}\right] \leq C_T (n!)^{  \beta(2-2\kappa_1)} C^n     T^{ 2\beta \kappa_1 n }
\end{aligned}
\end{equation*}
From this bound we conclude  that there exists a constant
$C_{d,T,k}>0$ such that
\begin{eqnarray*}
&&\mathbb{E}\left[\exp\left\{C_{d,T,k}\left|\widehat{\alpha}^{(k)}(0)\right|^{\beta}\right\}\right]
 =\sum_{n=0}^{\infty}\frac{C_{d,T,k}^n}{n!}\EE\left|\widehat{\alpha}^{(k)}(0)
\right|^{n\beta}\\
&&\qquad\le C_T \sum_{n=0}^\infty C_{d,T, k}^n (n!)^{  \beta(2-2\kappa_1)-1} C^n     T^{ 2\beta \kappa_1 n }   <\infty,\\
\end{eqnarray*}
when $C_{d, T, k}$ is sufficiently small (but strictly positive),
where  $\beta=\frac{H_1+H_2}{2dH_1H_2}$.

%According to Theorem 2, $\mathbb{E}[(\hat{\alpha}^{(k)}(0))^2]$ is finite  if
%$\frac{1}{H_1d}+\frac{1}{H_2d}>1+2k$. In particular,  the condition of
% $\hat{\alpha}^{(2)}(0)\in L^2(\Omega)$ is $\frac{1}{H_1d}+\frac{1}{H_2d}>5$, which
% is the same condition in Theorem 1. Comparing it with the
%condition of $I(B^H,\tilde{B}^H)$ existence in $L^2$ in Nualart et al. (2007), where
%it need $Hd<2$ for two independent fractional Brownian motions taking the same  Hurst index
%H, our condition  also covers the result
%when only take $k=0$.

%Lastly, we discuss some integrability conditions  of derivative local time of fractional
% Brownian motions. For making comparisons, we find  that there is some difference between
%  difference order of derivative  local time.
%
%\begin{theorem}
%Let $B^{H_1}$ and $\widetilde{B}^{H_2}$  be two independent d-dimensional fractional
% Brownian motions. Then,
%\begin{equation*}
%\begin{aligned}
%\mathbb{E}\left|\hat{\alpha}^{(k)}(0)\right|<\infty,
%\end{aligned}
%\end{equation*}
% if and only if Hurst parameters $H_1$  and $H_2$ satisfy
% \begin{equation}
% \frac{1}{H_1 }+\frac{1}{H_2 }> d+|k|\label{e.cond_H}
% \end{equation}
%\end{theorem}

%\begin{proof}
%From Theorem 1, we see that the condition \eqref{e.cond_H}
%is sufficient. Now we show that it is also necessary.

\medskip
\noindent{\bf Proof of part (iii)}. \
Without loss generality, we consider only the case
$k=(k_1,0, \cdots, 0)$ and we denote $k_i$ by $k$.
By the definition of k-order derivative local time of independent d-dimensional fractional
 Brownian motions, we have
\begin{equation*}
\begin{aligned}
\mathbb{E}\left[\hat{\alpha}^{(k)}_{\varepsilon}(0)\right]
&=\frac{1}{(2\pi)^d}\int_{0}^T\int_0^T\int_{\mathbb{R}^d}\EE \left[e^{i\langle
\xi,B^{H_1}_t-\widetilde{B}^{H_2}_s\rangle}\right]
e^{-\frac{\varepsilon\mid\xi\mid^2}{2}}\mid\xi_1\mid^kd\xi dtds\\
&=\frac{1}{(2\pi)^d}\int_{0}^T\int_0^T\int_{\mathbb{R}^d}e^{-(\varepsilon+t^{2H_1}
+s^{2H_2})
\frac{\mid\xi\mid^2}{2}}
\mid\xi_1\mid^kd\xi dtds.\\
\end{aligned}
\end{equation*}
Thus, we have
\begin{equation*}
\mathbb{E} \left[\hat{\alpha}^{(k)}(0)\right]
=\frac{1}{(2\pi)^d}\int_{0}^T\int_0^T\int_{\mathbb{R}^d}e^{-(t^{2H_1}+s^{2H_2})
\frac{\mid\xi\mid^2}{2}}
\mid\xi_1\mid^kd\xi dtds.
\end{equation*}
%According to the fact
% \begin{equation*}
%\int_{\mathbb{R}}x^ke^{-\frac{x^2(t^{2H_1}+s^{2H_2})}{2}}dx=
%  \left\{
%   \begin{aligned}
%   &\sqrt{2\pi}(k-1)!!(t^{2H_1}+s^{2H_2})^{-\frac{k+1}{2}}, if~~ k ~~is ~~even, \\
%   &0, ~~~~~~~~~~~~~~~~~~~~~~~~~~~~~~~~~~~~~~~if ~~k ~~is ~~odd,\\
%   \end{aligned}
%   \right.
%  \end{equation*}
Integrating with respect to $\xi$, we find
\begin{equation*}
\mathbb{E}\left[\hat{\alpha}^{(k)}(0)\right]
=c_{k, d} \int_{0}^T\int_0^T(t^{2H_1}+s^{2H_2})
^{-\frac{(k+d) }{2}}dtds
\end{equation*}
for some constant $c_{k, d}\in (0,\infty)$.

We are going to deal with  the above integral.  Assume first $0<H_1\le H_2<1$.  Making substitution  $t=u^{\frac{H_2}{H_1}}$
yields
\begin{eqnarray}
I_4
&:=&\int_{0}^T\int_0^T(t^{2H_1}+s^{2H_2})
^{-\frac{(k+d) }{2}}dtds \nonumber\\
&=& \int_{0}^T\int_0^{T^{\frac{H_1}{H_2}}}(u^{2H_2}+s^{2H_2})^{-\frac{k+d}{2}}
u^{\frac{H_2}{H_1}-1}duds\,.
\end{eqnarray}
Using  polar   coordinate
$u=r\cos\theta$ and $s=r\sin\theta$, where $0\leq\theta\leq\frac{\pi}{2}$ and
$0\leq r\leq T$  we have
\begin{eqnarray}
I_4
&\ge &  \int_0^{\frac{\pi}{2}}(\cos\theta)^{\frac{H_2}{H_1}-1}
(\cos^{2H_2}\theta+\sin^{2H_2}\theta)^{-\frac{(k+d) }{2}} d\theta\int_0^{ T
^{\frac{H_1}{H_2}}}
r^{-(k+d)H_2 +
\frac{H_2}{H_1}}dr  \nonumber\\
\label{e.11}
\end{eqnarray}
since the planar  domain $\left\{(r, \theta), 0\le r\le
T\wedge T^{\frac{H_1}{H_2}}, 0\le \theta\le \frac{\pi}{2}\right\}$ is contained in
the planar domain $\left\{(s,u), 0\le s\le T\,, 0\le u\le T
^{\frac{H_1}{H_2}} \right\}$.
The integral with respect to $r$  appearing in  \eqref{e.11}
 is finite  only if  $-(k+d)H_2 +\frac{H_2}{H_1}>-1$,   namely,
only when the condition \eqref{e.cond_H} is satisfied.  The case
$0<H_2\le H_1<1$ can be  dealt  similarly.  This completes the proof of
our main  theorem. \hfill $\Box$

%
%
%
% if and only if Hurst parameters $H_1$  and $H_2$ satisfy $\frac{1}{H_1d}+\frac{1}{H_2d}>1+k$. )
%
%Secondly, if $0<H_2<H_1<1$, taking $s={u}^{\frac{H_1}{H_2}}$, there is
%\begin{equation}
%\begin{aligned}
%&\int_{0}^T\int_0^T(t^{2H_1}+s^{2H_2})
%^{-\frac{(k+1)d}{2}}dtds\\
%&~~=\int_{0}^{T^{\frac{H_2}{H_1}}}\int_0^{T}(t^{2H_1}+u^{2H_1})^{-\frac{(k+1)d}{2}}
%u^{\frac{H_1}{H_2}-1}dtdu\\
%&~~\leq\int_0^{\frac{\pi}{2}}(\cos\theta)^{\frac{H_1}{H_2}-1}
%(\cos^{2H_1}\theta+\sin^{2H_1}\theta)^{-\frac{(k+1)d}{2}}d\theta\int_0^{\sqrt{2}T}
%r^{-(k+1)H_1d+
%\frac{H_1}{H_2}}dr\\
%&~~<\infty,\\
%\end{aligned}
%\end{equation}
%if and only if $-(k+1)H_1d+\frac{H_1}{H_2}>-1$, i.e., $\frac{1}{H_1d}+\frac{1}{H_2d}>1+k$.
%
%Lastly, when $H_1=H_2=H$, we have
%\begin{equation}
%\begin{aligned}
%&\int_{0}^T\int_0^T(t^{2H}+s^{2H})
%^{-\frac{(k+1)d}{2}}dtds\\
%&~~\leq\int_0^{\frac{\pi}{2}}
%(\cos^{2H}\theta+\sin^{2H}\theta)^{-\frac{(k+1)d}{2}}d\theta\int_0^{\sqrt{2}T}
%r^{1-(k+1)dH}dr\\
%&~~<\infty,\\
%\end{aligned}
%\end{equation}
%if and only if $1-(k+1)Hd>-1$, in other word, it need $Hd<\frac{2}{k+1}$.
%
%Therefore,
%\begin{equation*}
%\begin{aligned}
%\mathbb{E}[\hat{\alpha}^{(k)}(0)]<\infty,
%\end{aligned}
%\end{equation*}
% if and only if  $\frac{1}{H_1d}+\frac{1}{H_2d}>1+k$.
%\end{proof}
% In Proposition 1, we obtain sufficient and necessary condition of  $\mathbb{E}
% [\hat{\alpha}^{(k)}(0)]<\infty$, which is different from the condition in the
%  proof of Theorem 2 for $n=1$.


\begin{thebibliography}{10}
%\bibitem{[1]} C. Bender. An It\^{o} formula for generalized of
%a fractional Brownian motion with arbitrary Hurst parameter.
%Stoch. Proce. Their Appl., 2003, 104: 81-106.
%\bibitem{[2]} F. Biagini, Y. Hu, B. $\Phi$ksendal
%,et al.
%Stochastic calculus for fractional Brownian motion and applications.
%London: Springer-Verlag, 2008.
\bibitem{nondet} S. M.  Berman. Local nondeterminism and local times of Gaussian processes. Indiana Univ. Math. J., 1973, 23: 69-94.
\bibitem{jung1} P. Jung, G. Markowsky. On the Tanaka formula for the derivative of
self-intersection local time of fractional Brownian motion. Stoc. Proces. Their
Appl., 2014, 124: 3846-3868.
\bibitem{jung2} P. Jung, G. Markowsky. H$\ddot{o}$lder continuity and occupation-time
formulas for fBm
self-intersection local time and its derivative. J. Theor. Probab., 2015, 28: 299-312.
\bibitem{hu17} Y. Hu. Analysis on Gaussian spaces. World Scientific, 2016.
\bibitem{hunualart} Y. Hu, D. Nualart. Renormalized self-intersection local time for
 fractional Brownian motion. Ann. Pro., 2005, 33(3): 948-983.
\bibitem{hunualartsong} Y. Hu, D. Nualart, J. Song. Integral representation of renormalized
self-intersection local times. J. Func. Anal., 2008, 255(9): 2507-2532.
\bibitem{HHNT}Y. Hu, J. Huang, D. Nualart, S. Tindel. Stochastic heat equations with general multiplicative Gaussian noises: H\"older continuity and intermittency. Electron. J. Probab., 2015,  20(55): 1-50.
\bibitem{nualartortiz} D. Nualart, S. Ortiz-Latorre.
Intersection local time for two  independent fractional Brownian
motions. J. Theor. Prob., 2007, 20: 759-767.
\bibitem{oliveira} M. Oliveira, J. Silva, L. Streit. Intersection local
times of independent fractional Brownian motions as generalized
white noise functionals. Acta Appl. Math., 2011, 113(1): 17-39.
\bibitem{wuxiao} D. Wu, Y. Xiao. Regularity of intersection local times of fractional
Brownian motions. J. Theor. Prob., 2010, 23(4): 972-1001.
\bibitem{yan} L. Yan. Derivative for the intersection local time of fractional Brownian
motions. arXiv:1403.4102v3, 2014.
\bibitem{yanyu} L. Yan, X. Yu. Derivative for self-intersection local time of
multidimensional fractional Brownian motion. Stoch., 2015, 87(6): 966-999.
\end{thebibliography}
\end{document}